\documentclass[reqno]{amsart}
\usepackage{amssymb}
\usepackage{amsmath} 
\usepackage{epsf} 
\usepackage{latexsym} 
\usepackage{amscd}
\usepackage{graphics}
\usepackage{graphpap}

\newtheorem{thm}{Theorem}[section]
\newcommand{\bthm}{\begin{thm}}
\newcommand{\ethm}{\end{thm}}

\newtheorem{prop}[thm]{Proposition}
\newcommand{\bprp}{\begin{prop}}
\newcommand{\eprp}{\end{prop}}

\newtheorem{fact}[thm]{Fact}
\newcommand{\bfct}{\begin{fact}}
\newcommand{\efct}{\end{fact}}

\newtheorem{prob}[thm]{Problem}
\newcommand{\bprb}{\begin{prob}}
\newcommand{\eprb}{\end{prob}}

\newtheorem{quest}[thm]{Question}
\newcommand{\bqtn}{\begin{quest}}
\newcommand{\eqtn}{\end{quest}}

\newtheorem{lem}[thm]{Lemma}
\newcommand{\blem}{\begin{lem}}
\newcommand{\elem}{\end{lem}}

\newtheorem{claim}[thm]{Claim}
\newcommand{\bclm}{\begin{claim}}
\newcommand{\eclm}{\end{claim}}

\newtheorem{cor}[thm]{Corollary}
\newcommand{\bcor}{\begin{cor}}
\newcommand{\ecor}{\end{cor}}

\newtheorem{conj}[thm]{Conjecture}
\newcommand{\bcnj}{\begin{conj}}
\newcommand{\ecnj}{\end{conj}}

\theoremstyle{definition}
\newtheorem{defn}[thm]{Definition}
\newcommand{\bdfn}{\begin{defn}}
\newcommand{\edfn}{\end{defn}}

\newtheorem{spec}[thm]{Specializing}
\newcommand{\bspc}{\begin{spec}}
\newcommand{\espc}{\end{spec}}

\theoremstyle{remark}
\newtheorem{rem}[thm]{Remark}
\newcommand{\brem}{\begin{rem}}
\newcommand{\erem}{\end{rem}}

\newtheorem{cnv}[thm]{Convention}
\newcommand{\bcnv}{\begin{cnv}}
\newcommand{\ecnv}{\end{cnv}}

\newtheorem{exam}[thm]{Example}
\newcommand{\bexm}{\begin{exam}}
\newcommand{\eexm}{\end{exam}}

\newtheorem{exercise}[thm]{Exercise}
\newcommand{\bexr}{\begin{exercise}}
\newcommand{\eexr}{\end{exercise}}

\def\AA {{\mathbb A}}

\def\NN {{\mathbb N}}

\def\QQ {{\mathbb Q}}
\def\RR {{\mathbb R}}
\def\TT {{\mathbb T}}

\def\ZZ {{\mathbb Z}}

\def\sB {{\mathcal B}}

\def\implies {\Longrightarrow}

\def\iff {\Longleftrightarrow}
\def\and {\wedge}
\def\qed {{\blacksquare}}
\def\onto {\,\rule[.04in]{.15in}{.01in}\kern-4pt\raise.3ex\hbox{$\scriptscriptstyle>\!>\,$}}
\def\to {\longrightarrow}

\def\lc {locally compact}
\def\st {such that}
\def\tg {topological group}
\def\sii {if and only if}

\def\nhd {{neighborhood}}

\def\ie {{\em i.e.,}}

\def\< {{\langle}}
\def\> {{\rangle}}

\def\gs {{\sigma}}

\def\pf {{\em \noindent Demostraci\'on:}}
\def\epf {~\hfill $\qed$}

\def\pf {{\em \noindent Proof:}} 
\def\l( {{\left)}} 
\def\r( {{\right)}} 
\def\l[ {{\left[}} 
\def\r] {{\right]}} 
\def\l{ {{\left{}} 
\def\r} {{\right}}} 

\def\mkp{{\medskip}}

\begin{document}

\title[When a totally bounded group topology is a Bohr Topology]{When a totally bounded group topology is the Bohr Topology of a LCA group}
\date{\today}
\author{Salvador Hern\'andez}
\address{Departmento of Matem\'aticas\\
Universitat Jaume I, 8029-AP\\
Castell\'on, Spain}
\email{hernande@mat.uji.es}
\author{F. Javier Trigos-Arrieta}
\address{Department of Mathematics\\
California State University, Bakersfield\\
Bakersfield, CA 93311}
\email{jtrigos@csub.edu}

\thanks{ The first listed author acknowledges partially support by the Spanish Ministerio de Econom\'{i}a y Competitividad,
grant MTM2016-77143-P (AEI/FEDER, EU),
and the Universitat Jaume I, grant P1171B2015-77.}


\begin{abstract}
We look at the Bohr topology of maximally almost periodic groups (MAP, for short). Among other results, we investigate when
a totally bounded abelian group $(G,w)$ is the Bohr reflection of a locally compact abelian group. Necessary and sufficient
conditions are established in terms of the inner properties of $w$. As an application, an example of a MAP group $(G,t)$  is given
such that every closed, metrizable subgroup $N$ of $bG$ with $N \cap G = \{0\}$ preserves compactness but $(G,t)$ does not strongly
respects compactness. Thereby, we respond to Questions 4.1 and 4.3 in \cite{comftrigwu}.
\end{abstract}

\maketitle

\section{Introduction}

For each topological group $(G,t)$ there is associated a compact Hausdorff group
$bG$ and a continuous homomorphism $b$ from $(G,t)$ onto a dense subgroup of $bG$
with the following universal property:\ for every continuous homomorphism $h$
from $(G,t)$ into a compact group $K$ there is a continuous homomorphism $h^+$
from $bG$ into $K$ such that $h=h^{+}\circ b$. The group $bG$ is essentially
unique;\ it is called the Bohr compactification of $G$ (see Heyer \cite{Heyer1970}
for a careful examination of $bG$ and its properties). 
Here, we restrict our attention to those
groups such that the homomorphism $b$ above is one-to-one; these are exactly
the {\em maximally almost periodic} (\textit{MAP}) groups.

For such a topological group $(G,t)$ we
denote by $(G,t^{+})$ the underlying group of $G$ equipped with its Bohr
topology. Evidently, $(G,t^{+})$ is an example of a totally bounded group.
The following notion plays an important r\^{o}le for the rest of this
discussion: A MAP group $(G,t)$ is said to \textit{respect} a topological
property $\mathcal{P}$ if a subset $A$ of $G$ has $\mathcal{P}$ as a
subspace of $(G,t)$ if and only if $A$ has $\mathcal{P}$ as a subspace of $(G,t^{+})$
(Trigos-Arrieta \cite{Trigos-Arrieta1991}).

The question of the disposition or placement of a LCA group $(G,\tau)$ within its
Bohr compactification $bG$ has been investigated by many workers. It is
known for such $G$, for example, that $G^{+}$ is \textit{sequentially closed}
in $bG$ in the sense that no sequence from $bG$ can converge to a point of $%
bG\backslash G$ \cite{reid}. And Glicksberg \cite{glicks1962} has shown that
LCA groups respect compactness. This result concerning LCA
groups is one of the pivotal results of the subject, often referred to as
\emph{Glicksberg's theorem}. Hughes \cite{Hughes1973} proved a generalization of
Glicksberg's theorem to (not necesary Abelian) locally compact groups by
considering the weak topology generated by the continuous irreducible
unitary group representations. Several authors have achieved additional
results which continue the lines of investigation suggested above (see \cite{GalHerWu}).
Nevertheless many questions relating the topology of a MAP group
with its Bohr topology are still open in general. In this paper,
we continue with the investigation of the Bohr topology of maximally almost periodic groups groups accomplished in \cite{comftrigwu} and \cite{Galindo2004}.
Among other results, we investigate when
a totally bounded abelian group $(G,w)$ is the Bohr reflection of a locally compact abelian group. Necessary and sufficient
conditions are established in terms of the inner properties of $w$. As an application, an example of a MAP group $(G,t)$ is given
such that every closed, metrizable subgroup $N$ of $bG$ with $N \cap G = \{0\}$ preserves compactness but $(G,t)$ does not strongly
respects compactness (see Definition \ref{str_resp_cpt} below). Thereby, we respond to Questions 4.1 and 4.3 in \cite{comftrigwu}.

\section{Preliminaries}

If $G$ is an Abelian group, the symbol $w$ stands for a Hausdorff precompact group topology on $G$, and $\tau$ stands for a Hausdorff locally compact group topology on $G$. Then $\tau^+$ stands for the Bohr  topology on $G$, associated to $\tau$. The topology $\tau^+$ is a Hausdorff precompact group topology on $G$, weaker than $\tau$. If $(G,t)$ is a \tg \ and $H$ is a subgroup of $G$, then $(H,t)$ stands for the \tg \ $H$ equipped with the inherited topology from $(G,t)$. $(H,t)<(G,t)$ means that $(H,t)$ is a closed subgroup of $(G,t)$. If so, $(G/H,t)$ stands for the natural quotient group. Let $\sB:=\{(G,w):\exists \tau [w=\tau^+]\}$.  Members of the class $\sB$ will be called {\em Bohr groups}. If $\tau$ is the discrete topology, we write $\#$ instead of $\tau^+$. We identify the torus $\TT$ with the group $[0,1) \subset \RR$ equipped with the operation $+$~mod~$1$. If $(G,t)$ is a \tg , the symbol $(G,t)\widehat{\:\:}$ stands for the $t$-continuous homomorphisms from $G$ to $\TT$. $(G,t)\widehat{\:\:}$ is a group, and it is called {\em the group of characters of $(G,t)$}. Notice that $(G,w) \in \sB \implies (G,w)\widehat{\:\:}=(G,\tau^+)\widehat{\:\:}=(G,\tau)\widehat{\:\:}$, for some \lc \ \tg \ topology $\tau$. By $((G,t)\widehat{\:\:},w)$ we mean the group $(G,t)\widehat{\:\:}$ equipped with the finite-open topology. Similarly, by $((G,t)\widehat{\:\:},k)$ we mean the group $(G,t)\widehat{\:\:}$ equipped with the compact-open topology. If $H$ is a subgroup of $G$, {\em the annihilator of $H$ in $(G,t)\widehat{\:\:}$}, denoted by $\AA((G,t)\widehat{\:\:},H)$ is the {\em subgroup of} $(G,t)\widehat{\:\:}$ consisting of those characters $\varphi$ from $G$ to $\TT$, \st \ $\varphi[H]=\{0\}$. If $X$ is a group of homomorphisms from $G$ to $\TT$, we say that $X$ {\em separates the points of $G$} if whenever $g \in G$ and $g \neq 0_G$, the identity of $G$, then there is $\varphi \in X$ \st \ $\varphi(g) \neq 0$. Hom$(G,\TT)$ stands for the group of all group homomorphisms from $G$ to $\TT$.

\begin{thm}\label{prelims}
\begin{enumerate}
\item $((G,\tau)\widehat{\:\:},k)$ {\em is \lc . $\tau$ discrete (compact, resp.) implies $((G,\tau)\widehat{\:\:},k)$ compact (discrete, resp.).
\item $\Omega:(G,\tau) \to (((G,\tau)\widehat{\:\:},k)\widehat{\:\:},k), g\mapsto (\Omega(g):((G,\tau)\widehat{\:\:},k)\to \TT, \phi \mapsto \phi(g))$ is a topological isomorphism onto. $\Omega$ is called {\em the evaluation map (of $(G,\tau)$).}
\item The group $($}Hom$(G,\TT),w)$ {\em is compact.
\item If $X \subseteq$}~Hom$(G,\TT)${\em, then $(X,w)$ is precompact.
\item If $X = (G,w)\widehat{\:\:}$ separates the points of $G$ and $t$ is the weakest topology on $G$ that makes the elements of $X$ continuous, then $(G,w)=(G,t)$.
\item If $X \subseteq$}~Hom$(G,\TT)$ {\em separates the points of $G$, then $(X,w)$ is a dense subgroup of} $($Hom$(G,\TT),w)$.
\item $\omega:(G,w) \to (((G,w)\widehat{\:\:},w)\widehat{\:\:},w), g\mapsto (\omega(g):((G,w)\widehat{\:\:},w)\to \TT, \phi \mapsto \phi(g))$ {\em is a topological isomorphism onto. $\omega$ is called}  the evaluation map (of $(G,w)$).
\item {\em If $\tau$ is \lc , then $K$ is a compact subspace of $(G,\tau) \iff K$ is a compact subspace of $(G,\tau^+)$.}
\item {\em If $\tau$ is \lc \ and $(H,\tau)<(G,\tau)$, then the Bohr topology of $(H,\tau)$ as a \lc \ \tg \ is the same as the topology that it inherits as a subgroup of $(G,\tau^+)$.}
\end{enumerate}
\end{thm} 

\pf \ 1 and 2 are the celebrated Pontryagin-van Kampen theorem. 3 follows from 1 taking $\tau$ discrete. 4 and 5 follow from Comfort and Ross \cite{ComfortRoss1964}. 6 is \cite{hewitt1963abstract} (26.16). 7 is Raczkowski and Trigos-Arrieta \cite{racz-trig2001}. 8 is Glicksberg's Theorem \cite{glicks1962}. 9 is done in \cite{Trigos-Arrieta1991}. \epf

\section{Necessary conditions.}

\begin{thm} \label{G_in_Gbar_densely} If $(G,w)$ is \st \ $((G,w)\widehat{\:\:},k)$ is \lc , then $(\overline{G}, \tau):=(((G,w)\widehat{\:\:},k)\widehat{\:\:},k)$ satisfies that $(G,w)$ is contained in $(\overline{G}, \tau^+)$ densely.
\end{thm} 

\pf \ The identity $I:((G,w)\widehat{\:\:},k)\to((G,w)\widehat{\:\:},w)$, is clearly continuous and onto. Hence, the adjoint map, which is the containment, $\widehat{I}:(((G,w)\widehat{\:\:},w)\widehat{\:\:},w) \to (((G,w)\widehat{\:\:},k)\widehat{\:\:},w)$ is continuous. By Theorem \ref{prelims}.7, $(((G,w)\widehat{\:\:},w)\widehat{\:\:},w)=(G,w)$. By Theorem \ref{prelims}.1, the group $(\overline{G}, \tau)$ is \lc . Hence, $(\overline{G}, \tau^+)=(((G,w)\widehat{\:\:},k)\widehat{\:\:},w)$. That the containment is dense, follows from Theorem \ref{prelims}.6.\epf

\blem\label{Gw_inB=>Gw^k_loc_cpt} If $(G,w) \in \sB$, then $((G,w)\widehat{\:\:},k)$ is locally compact.
\elem 

\pf \ Suppose that $\tau$ is \lc \ with $\tau^+=w$. By Theorem \ref{prelims}.8, $((G,w)\widehat{\:\:},k)=((G,\tau)\widehat{\:\:},k)$, which is \lc \ by Theorem \ref{prelims}.1. \epf

\bthm \label{HinB_G/HinB=>GinB} If $(H,w)<(G,w)$, $(H,w) \in \sB$, and $(G/H,w) \in \sB$, then $(G,w) \in \sB$.
\ethm 

\pf \ By Lemma \ref{Gw_inB=>Gw^k_loc_cpt}, both $((G/H,w)\widehat{\:\:},k)$ and $((H,w)\widehat{\:\:},k)$ are \lc . Let $X:=((G,w)\widehat{\:\:},k)$. We have that $(G/H,w)\widehat{\:\:}=\AA(X,H)$ \cite{hewitt1963abstract} (23.25 \& 23.30) and $(H,w)\widehat{\:\:}=X/\AA(X,H)$ \cite{hewitt1963abstract} (24.11 \& 23.30). By \cite{hewitt1963abstract} (5.25), we have that $X$ is \lc , hence $(\overline{G},\tau):= (X\widehat{\:\:},k)$ is \lc \ as well. It follows that $G\subseteq (\overline{G},\tau)$ densely, since both $G$ and $\overline{G}$ separate the points of $X$. We claim that $H$ as a subgroup of $(\overline{G}, \tau)$ is closed. For, $(H,w) \in \sB \implies \exists \tau_H$ \lc \ \st \ $(H,\tau_H^+)=(H,w)$. The latter implies that $(H,\tau_H)\widehat{\:\:}=(H,\tau_H^+)\widehat{\:\:}=(H,w)\widehat{\:\:}=X/\AA(X,H)$. Therefore, $(H,\tau_H)=(X/\AA(X,H))\widehat{\:\:} < (\overline{G},\tau)$ which proves the claim. Similarly, we claim that $G/H$ is \lc \ as a subgroup of $(\overline{G}, \tau)/H$. For, $(G/H,w) \in \sB \implies \exists \tau_{G/H}$ \lc \ \st \ $(G/H,\tau_{G/H}^+)=(G/H,w)$. The latter implies that $(G/H,\tau_{G/H})\widehat{\:\:}=(G/H,\tau_{G/H}^+)\widehat{\:\:}=(G/H,w)\widehat{\:\:}=\AA(X,H)$. Therefore, $(G/H,\tau_{G/H})=\AA(X,H)\widehat{\:\:}=\overline{G}/H$; proving the second claim. Since $H$ and $G/H$ are \lc \ as a subgroup and quotient of $G<(\overline{G}, \tau)$ respectively, it follows, by \cite{hewitt1963abstract} (5.25), that $G$ is \lc \ as a subgroup of $(\overline{G}, \tau)$. Because $G$ is dense in $\overline{G}$, it follows that $G=\overline{G}$ \cite{hewitt1963abstract} (5.11), hence, $(G,\tau^+)=(G,w)$, as required. \epf

\bthm \label{cpt=finite} Assume all compact subsets of $(G,w)$ are finite. Then TFAE:
\begin{enumerate}
\item $(G,w) \in \sB$
\item $w=\#$
\item $((G,w)\widehat{\:\:},w)$ is compact.
\item $((G,w)\widehat{\:\:},k)$ is compact.
\item $((G,w)\widehat{\:\:},w)$ is \lc .
\item $((G,w)\widehat{\:\:},k)$ is \lc .
\item Every homomorphism $f:(G,w) \to \TT$ is continuous.
\end{enumerate}
\ethm 

\pf \ The hypothesis on the compact subsets of $(G,w)$ implies that $((G,w)\widehat{\:\:},w)=((G,w)\widehat{\:\:},k)$. Therefore $3 \iff 4$ and $5 \iff 6$. $2 \implies 1$ is obvious, Lemma \ref{Gw_inB=>Gw^k_loc_cpt} yields $1 \implies 6$. Because every group of characters equipped with the finite-open topology is precompact, $5 \implies 3$. Obviously, $2 \iff 7$. $3 \implies 7$ can be seen by using (3), (4) and (6) of Theorem \ref{prelims}.\epf

By \cite{hewitt1963abstract} (5.14), every \lc \ Abelian group $(G,\tau)$ contains an open compactly generated subgroup $(H,\tau)$. Notice then that $(G/H,\tau)$ is discrete. By properties of the Bohr topology, we have then that $(G,\tau^+)$ contains a compactly generated subgroup $(H,\tau^+)$ such that $(G/H,\tau^+)=(G/H,\#)$. This proves the sufficiency of the following:

\blem \label{Gw_in_B_sii_compactly_generated} $(G,w) \in \sB$ \sii \ there is a compactly generated subgroup $(H,w)$ of $(G,w)$, \st \ $(H,w) \in \sB$ and  $(G/H,w)= (G/H,\#)$.
\elem 

\pf \ For the necessity, use Theorem \ref{HinB_G/HinB=>GinB}. \epf

Obviously:

\bcor \label{Hwnot_in_B_or_G/Hwnot_in_B=>G_not_in_B} If for every compactly generated subgroup $(H,w)$ of $(G,w)$ we have that $(H,w)\not \in \sB$ or $(G/H,w)\neq (G/H,\#)$,
then $(G,w) \not\in \sB$. \epf
\ecor 

By \cite{hewitt1963abstract} (9.8) a compactly generated \lc \ Abelian group must have the form $\ZZ^m\times \RR^n \times K$ where $m,n$ are non-negative integers and $K$ is a compact Abelian group. Because of Theorem \ref{prelims}.8 and properties of the Bohr topology, a compactly generated precompact group $(G,w) \in \sB$ \sii \ it has the form $(\ZZ^\#)^m\times (\RR^+)^n \times K$ where $m,n$ are non-negative integers and $K$ is a compact Abelian group, and by duality properties, if $(G,w) \in \sB$, then $((G,w)\widehat{\:\:},k)$ has the form $\TT^m\times \RR^n \times D$, where $m,n$ are non-negative integers and $D$ is a discrete Abelian group.

\bthm \label{Gw_compactly_generated=>GwinB_if...} Assume that $(G,w)$ is compactly generated. Then $1 \iff 2 \implies 3$:
\begin{enumerate}
\item $(G,w) \in \sB$,
\item $(G,w) =(\ZZ^\#)^m\times (\RR^+)^n \times K$, where $m,n$ are non-negative integers and $K$ is a compact Abelian group,
\item $((G,w)\widehat{\:\:},k)=\TT^m\times \RR^n \times D$ where $m,n$ are non-negative integers and $D$ is a discrete Abelian group.\epf
\end{enumerate}
\ethm 

(3) does not imply (2) above: Setting $(G,w):=\QQ^+$, one has $((G,w)\widehat{\:\:},k)=\RR$.

\bthm \label{Gw_cpt_gen+0dim=>Gw_inBsii} Assume that $(G,w)$ is compactly generated and 0-dimensional. Then $(G,w) \in \sB$ \sii \ the subgroup $B$ of $(G,w)$ of all compact elements of $(G,w)$ \cite{hewitt1963abstract} (9.10) is compact, and $(G/B,w)=(G/B,\#)$.
\ethm 

\pf \ $(\implies)$ Assume that $(G,w)=(G,\tau^+)$. By \cite{hewitt1963abstract} (9.8) there are $m,n\in \omega$ and a compact group $K$ \st \ $(G,\tau)=\RR^m\times \ZZ^n\times K$. By hypothesis, $m=0$ and $K$ is 0-dimensional. It follows that $(G,w)=(G,\tau^+)=(\ZZ^\#)^n\times K$,  $B=K$, and $(G/B,w)=(\ZZ^\#)^n$, as required. For $(\Longleftarrow)$, apply Theorem \ref{HinB_G/HinB=>GinB}. \epf

\blem \label{prod_in_B_sii_factors_in_B}
Suppose the groups $A, B$ and $G$ satisfy $(G,w)=(A,w) \times (B,w)$. Then $(G,w) \in \sB$ \sii \ both $(A,w) \in \sB$ and $(B,w)\in \sB$.
\elem 

\pf \ By the main hypothesis, $((G,w)\widehat{\:\:},k)=((A,w)\widehat{\:\:},k) \times ((B,w)\widehat{\:\:},k)$. If $(G,w) \in \sB$, say $(G,w) =(G,\tau^+)$, then, by Theorem \ref{prelims}.2, $(G,\tau)=(((G,w)\widehat{\:\:},k)\widehat{\:\:},k)=(((A,w)\widehat{\:\:},k)\widehat{\:\:},k) \times (((B,w)\widehat{\:\:},k)\widehat{\:\:},k)$. We write $(\overline{A},\tau):=(((A,w)\widehat{\:\:},k)\widehat{\:\:},k) $, and $(\overline{B},\tau):=(((B,w)\widehat{\:\:},k)\widehat{\:\:},k) $. Hence $(G,\tau)=(\overline{A},\tau) \times (\overline{B},\tau)$, and by Theorem \ref{G_in_Gbar_densely}, $A$ is a dense subgroup of $\overline{A}$, and $B$ is a dense subgroup of $\overline{B}$. By properties of the Bohr topology, $(A,w) \times (B,w)=(G,w) =(G,\tau^+)=(\overline{A},\tau^+) \times (\overline{B},\tau^+)$. This obviously implies that $(A,w) =(\overline{A},\tau^+)$ and $(B,w) =(\overline{B},\tau^+)$. The converse is obvious.\epf

\bthm \label{prod_in_B_iff_factors_nd_quoti_in_B} The \tg \ $(G,w) \in \sB$ \sii \ there are groups $A, B$ and $C$ \st \ (a) $G=A \times B$, (b) $(A\times\{0\},w)=(\RR^+)^n$ for some $n \in \omega$, and (c) $(\{0\}\times C,w)$ is a compact subgroup of $(\{0\}\times B,w)$ \st \ $(B/C,w)=(B/C,\#)$.
\ethm 

\pf \ $(\implies)$ Assume that $(G,w)=(G,\tau^+)$. By \cite{hewitt1963abstract} (24.30) there are $n\in \omega$, a \lc \ group $G_0$ and a compact group $K$ \st \ $(G,\tau)=\RR^n\times G_0$, and $(G_0/K,\tau)$ is discrete. Let $A, B$ and $C$ be the underlying groups of $\RR^n, G_0$, and $K$, respectively. By properties of the Bohr topology, we have that (a), (b) and (c) hold. For $(\Longleftarrow)$, apply Theorem \ref{HinB_G/HinB=>GinB} to see that $(\{0\}\times B,w) \in \sB$. That $(G,w) \in \sB$ follows after an application of Lemma \ref{prod_in_B_sii_factors_in_B}.\epf

\bthm \label{connect_comp} Consider the \tg \ $(G,w)$. Suppose that $F$ is its connected component. Then $(G,w) \in \sB$ \sii \ both $(F,w)\in \sB$, and $(G/F,w)\in \sB$.
\ethm 

\pf \ $(\implies)$ Assume that $(G,w)=(G,\tau^+)$. By \cite{hewitt1963abstract} (24.30) there are $n\in \omega$, a \lc \ group $G_0$ and a compact group $K$ \st \ $(G,\tau)=\RR^n\times G_0$, and $(G_0/K,\tau)$ is discrete. Let $C$ be the connected component of $G_0$. Clearly, $(C,\tau^+)<(G_0,\tau^+)=(G_0,w)$, with $(C,\tau^+)=(C,w)$ connected. Since $(G_0/C,\tau)$ is 0-dimensional, it follows, by properties of the Bohr topology \cite{Hernandez1998}, that $(G_0/C,\tau^+)$ is 0-dimensional. If $F:=(\RR^+)^n\times C$, it follows that $(F,w)\in \sB$ is the connected component of $(G,w)=(G,\tau^+)$. Of course $(G/F,w)=(G/F,\tau^+)\in \sB$. 
For $(\Longleftarrow)$, apply Theorem \ref{HinB_G/HinB=>GinB}.\epf

\bcor\label{Gw_in_B_TFAE} The following are equivalent, concerning a \tg \ $(G,w)$.
\begin{enumerate}
\item $(G,w) \in \sB$.
\item There are groups $A, B$ and $C$ \st \ (a) $G=A \times B$, (b) $(A\times\{0\},w)=(\RR^+)^n$ for some $n \in \omega$, and (c) $(\{0\}\times C,w)$ is a compact subgroup of $(\{0\}\times B,w)$ \st \ $(B/C,w)=(B/C,\#)$.
\item If $F$ is the connected component of $(G,w) \in \sB$ it follows that $(F,w)\in \sB$, and $(G/F,w)\in \sB$.
\item There is a compactly generated subgroup $(H,w)$ of $(G,w)$, \st \ $(H,w) \in \sB$ and  $(G/H,w)= (G/H,\#)$.
\end{enumerate}\ecor 

\bcor\label{cor_metrizable}
If $(G,w)$ is a metrizable totally bounded group that is in $\sB$, then $(G,w)$ is compact.
\ecor 


\section{Some topologies canonically associated to topological groups}

Following the terminology in \cite{comftrigwu}, for a (Hausdorff) space $X = (X,t)$
we denote 
by $kX$, or $(X, kt)$, the set $X$ with the topology $kt$ defined as follows:
A subset $U$ of $X$ is $kt$-open if and only if $U \cap K$ is (relatively) $t$-open in $K$ for every
$t$-compact subset $K$ of $X$. Then $kX$ is a $k$-space (that is, $kX = kkX$), $kt$
is the smallest $k$-space topology on $X$ containing $t$, 
it is the unique
$k$-space topology for $(X,t)$ such that $kt\supseteq t$, and the $kt$-compact sets
are exactly the $t$-compact sets. In like manner, we say that a map $f$ defined on $X$
is {\em $k$-continuous} when $f$ is continuous on each compact subset of $X$.\mkp

Given a topological abelian group $(G,t)$
with dual $X:=\widehat{(G,t)}$, for any subset $A$ of $G$, we
define $A^0:=\{\chi \in X: |\chi(g)|\leq 1/4 \:\forall \:g \in
A\}$. Assuming that we are considering the {\em dual pair}
$(G,X)$, for any subset $L$ of $X$, we define $L^0:=\{g\in G:
|\chi(g)|\leq 1/4 \:\forall \:\chi \in L\}$. This set operator
behaves in many aspects like the polar operator in vector spaces.
For instance, it is easily checked that $A^{000}=A^0$ for any
$A\subset G$. Given an arbitrary subset $A$ in $G$, we define the
\textit{quasi convex hull} of $A$, denoted $co(A)$, as the set
$A^{00}$. A set $A$ is said to be \textit{quasi convex} when it
coincides with its quasi convex hull. These definitions also apply
to subsets $L$ of $G$. The topological group $(G,t)$ is said to be
\textit{locally quasi convex} when there is a neighborhood base of
the identity consisting of quasi convex sets.

Let $(G,t)$ be a MAP \tg. In the sequel we are going to look at the following natural group topologies canonically attached to $(G,t)$.

\bdfn \label{defs}
\begin{enumerate}
\item {\bf The Bohr topology.} Denoted by $t^+$, it is the weak topology generated by the continuous homomorphisms from $G$ into $\TT$. It easily seen that the canonical map $b\colon (G,t)\mapsto (G,t^+)$ is an epireflective functor from the category of topological groups into the subcategory of totally bounded groups.

\item {\bf The locally quasi convex topology.} Denoted by $q[t]$, it is the finest locally quasi convex topology that is contained in $t$.
Again, it easily seen that the canonical map $q\colon (G,t)\mapsto (G,q[t])$ is an epireflective functor from the category of topological groups
into the subcategory of locally quasi convex groups.

\item {\bf The $g$-sequential topology.}  Denoted by $s_g[t]$, it is the finest group topology coarser than the {\em sequential modification of $t$,} \ie \ the largest topology on $G$ with the same $t$-convergent sequences.
When $t=s_g[t]$, it is said that $(G,t)$ is a \emph{$g$-sequential group}. In this case, the canonical map
$g\colon (G,t)\mapsto (G,s_g[t])$ defines a coreflective functor from the category of topological groups into the subcategory of $g$-groups.


\item {\bf The $k_g$-topology.}  Denoted by $k_g[t]$, it was originally defined by Noble \cite{noble70b} as the supreme of all group topologies on $G$ that lie between
$t$ and the $k$-topology attached to $t$. When $t=k_g[t]$, it is said that $(G,t)$ is a {\em $k_g$-group.} Also in this case, the canonical map $k_g\colon (G,t)\to (G,k_g[t])$ is a coreflective functor from the category of topological groups into the subcategory of $k_g$-groups.
\end{enumerate}
\edfn

First, we explore the relationship among these topologies. The topologies (3) and (4) have very similar properties and we will only present the
proofs for the $g$-sequential topology since basically the same proofs work for the $k_g$-topology.

\blem \label{lem_bn}
Let $(G,w)$ be a totally bounded group. Then $(G,w)$ is the Bohr reflection of a $g$-sequential group \sii \
$w$ is the finest among all totally bounded topologies on $G$ that share the same collection of convergent sequences.
If either condition holds, then $(G,s_g[w])$ has the same dual as $(G,w)$.
\elem 

\pf\ Assume that $w=\tau^+$ for some $g$-sequential topology $\tau$ on $G$. Let $\rho$ be another totally bounded topology on $G$
such that $w$ and $\rho$ have the same family of convergent sequences. Then the identity mapping $1_G\colon (G,s_g[w])\to (G,\rho)$
is sequentially continuous and, therefore, continuous as well. Now, since the map $g$ defines a correflective functor and $(G,\tau)$
is a $g$-sequential group, from $\tau\geq w$,
 it follows that $\tau\geq s_g[w]$, hence $1_G\colon (G,\tau)\to (G,\rho)$ is continuous.
Being $\rho$ totally bounded, this implies that $1_G\colon (G,\tau^+)\to (G,\rho)$ is continuous.
Since $\tau^+=w$, it follows that $w\geq \rho$.

Conversely, suppose that $w$ is the finest totally bounded topology on $G$ with the same family of convergent sequences.
First, we will see that $(G,s_g[w])$ has the same dual as $(G,w)$. Trivially, every $w$-continuous character is
automatically $s_g[w]$-continuous. We claim that if $\chi$ were a $s_g[w]$-continuous character that is not $w$-continuous, then $w\vee t_\chi$,
the supreme topology generated by $w$ and the initial topology generated by $\chi$, would be a totally bounded
topology with the same convergent sequences as $w$.

Indeed, let $(x_i)$ be a sequence in $G$ $w$-converging to some point $x_0\in G$. Since every $w$-convergent sequence is  $s_g[w]$-convergent,
it follows that $(x_{i})$ $s_g[w]$-converges to $x_0$. As a consequence $(\chi(x_{i}))$ converges to $\chi(x_0)$, and therefore,
$(x_{i})$ converges to $x_0$ in $w\vee t_\chi$. Thus  $w=w\vee t_\chi$,
which means that $\chi$ is $w$-continuous by Theorem \ref{prelims}.5. We have therefore verified that $w$ is the Bohr topology associated to $s_g[w]$.\epf

\blem \label{lem_bn2}
Let $(G,w)$ be a totally bounded group. Then $(G,w)$ is the Bohr reflection of a $k_g$-group \sii \
$w$ is the finest among all totally bounded topologies on $G$ that share the same collection of compact subsets.
If either condition holds, then $(G,k_g[w])$ has the same dual as $(G,w)$.
\elem

\bcor \label{cor_countable}
A countable totally bounded group $(G,w)$ is in $\sB$ \sii\ every character is continuous. Therefore, the groups $(G,s_g[w])$
and $(G,k_g[w])$ are discrete.
\ecor 

\pf \ ($\implies$) follows from Corollary \ref{Gw_in_B_TFAE}, while ($\Longleftarrow$) follows from Theorem \ref{prelims}.5 and Corollary \ref{Gw_in_B_TFAE}.2
That $s_g[w]$ and $k_g[w]$ are the discrete topology follows from Theorem \ref{cpt=finite}. \epf

\bcor \label{cor_Bohr1}
Let $(G,w)$ be a totally bounded group  that is the Bohr reflection of a $g$-sequential (resp. $k_g$) group $(G,\tau)$.
Then $w = s_g[w]^+$ (resp. $w = k_g[w]^+$).
\ecor 

\pf \ 
It suffices to notice that $(G,s_g[w])$ ($(G,k_g[w])$, resp.) has the same dual as $(G,w)$.\epf

\blem \label{lem_chacterization}
Let $(G,t)$ be a \tg \ and let $G':=(G,t)\widehat{\:\:}$ denote its dual group.
Then TFAE:
\begin{enumerate}
\item Every sequentially continuous character on $(G,t)$ is continuous;
\item  $(G,s_g[t])\widehat{\:\:}=G'$;
\item $(G,s_g[t]^+)=(G,t^+)$;
\item There exists a $g$-sequential topology $\tau$ on $G$ such that $t\subseteq \tau$ and $t^+=\tau^+$.
\end{enumerate} \elem 

\blem \label{lem_chacterization2}
Let $(G,t)$ be a \tg \ and let $G':=(G,t)\widehat{\:\:}$ denote its dual group.
Then TFAE:
\begin{enumerate}
\item Every $k$-continuous character on $(G,t)$ is continuous;
\item  $(G,k_g[t])\widehat{\:\:}=G'$;
\item $(G,k_g[t]^+)=(G,t^+)$;
\item There exists a $k_g$-topology $\tau$ on $G$ such that $t\subseteq \tau$ and $t^+=\tau^+$.
\end{enumerate} \elem

\bcor \label{cor_compactness}
The group $(G,t)$ respects convergent sequences \sii\ $t\subseteq s_g[t^+]$.
\ecor 

\bcor \label{cor_compactness2}
The group $(G,t)$ respects compact subsets \sii\ $t\subseteq k_g[t^+]$.
\ecor
\medskip

In order to characterize Bohr groups we need two basic notions. The first one is well known
and the later was introduced in \cite{HT_jmaa_2005}.\medskip

A family $\mathcal N$ of subsets of a topological space $X$ is a \emph{network at $x\in X$} if for every neighborhood $U$ of $x$
there exists an $M\in \mathcal N$ such that $x\in M\subseteq U$. If $\mathcal N$ is a network at each point in $X$,
we say that $\mathcal N$ is a \emph{network for $X$}. 
\medskip

For any topological group $(G,t)$, we say that $W \subset G$ is a \emph{$k$-neighborhood} of $0$
if for any $t$-compact subset $K \subset G$ containing $0$, there exists a neighborhood $U$ of $0$ such that
$U \cap K\subset W \cap K$. It is not true in general that if $x \in G$ and $U$ is a
$k$-\nhd \ of $x$ in $(G,t)$ then $U$ is a \nhd \ of $x$ in $(G,k_g[t]))$. However, when $U$ is a quasi
convex set the answer is positive (see Proposition 1 in \cite{HT_jmaa_2005}).\medskip

A topological space is said to be {\em hemicompact} if in the family of all the compact subspaces of $X$ ordered by $\subseteq$ there there is a countable cofinite subfamily. The concept was introduced by Arens in \cite{arens1946}. Hemicompact spaces are of course $\gs$-compact but $\QQ$ shows that the containment is proper. We now look at groups that are hemicompact. 


\bthm \label{th_characterization}
Let $(G,w)$ be a hemicompact, totally bounded group whose cardinality is not Ulam-measurable.
Then $(G,w)\in \mathcal B$ \sii \ the following properties hold:
\begin{enumerate}
\item  Every sequentially continuous character on $(G,w)$ is continuous.
\item There exists a compact subgroup $K$ of $G$ such that $G/K$ has a countable network at $0$ consisting of $k$-neighborhoods of $0$.
\end{enumerate}
\ethm 

\pf\ First, we notice that a wide use of duality techniques are essential for the proof. Assume that $(G,w)\in \mathcal B$ and let $\tau$ be a locally compact topology on $G$ such that $\tau^+=w$.
Then $(G,\tau)$ satisfies the two assertions above. Indeed, that $(G,\tau)$ satisfies (1) is due to results of Varopoulos \cite{varop64} and Reid \cite{reid}. 
On the other hand, the celebrated Kakutani-Kodaira Theorem \cite{comf-hbook} (3.7) and Theorem \ref{prelims}.8 imply that $G$ contains a compact subgroup $K$
such that $(G/K,\tau/K)$ is metric. Therefore $(G/K,\tau/K)$ will be hemicompact and metric. Now, every LCA group is a locally quasi convex $k_g$-group. Therefore,
there exists a countable neighborhood base at the neutral element $\mathcal N=\{U_n\}$ consisting of $k$-neighborhoods, quasi convex sets. Then $\mathcal N$ is a countable network for $G/K$,
which proves (2).

Conversely, suppose that (1) and (2) hold. By Theorem \ref{HinB_G/HinB=>GinB}, in order to prove that $(G,w)\in \mathcal B$, it will suffice to
verify that $(G/K,w/K)\in \mathcal B$. In other words, there is no loss of generality in assuming that $(G,w)$
has a countable network at $0_G$ consisting of $k$-neighborhoods of $0$.

Let us denote by $G':=((G,w)\widehat{\:\:},k),$ the dual topological group of $(G,w)$. Because $(G,w)$ is hemicompact, it follows that $G'$
is metric and, by (1), it follows that $G'$ is complete metric. On the other hand, by \cite[Lemma 5]{HT_jmaa_2005}, we have that if $F$ is a $k$-neighborhood
of $0_G$, then $F^0$ is precompact in $G'$. Furthermore, being $F^0$ closed in $G'$, which is complete, it follows that $F^0$ is in fact a compact subset of $G'$.
Since $(G,w)$ has a countable network $\{W_n\}$ at $0_G$ consisting of quasi convex subsets, it follows that $G'=\cup_{n\in\NN} W^0_n$.
Therefore $G'$ is $\sigma$-compact. Furthermore, since $G'$ is complete metric and $\sigma$-compact, by the Baire's category theorem, it follows that $G'$ is locally compact. { In particular, $G'$ is hemicompact \cite{arens1946}.}

{ If $\widehat{G'}$ denote the dual group of $G'$, then obviously $(\widehat{G'},k(\widehat{G'},G'))$ will be a locally compact abelian group, metrizable in addition since $G'$ is hemicompact. 
We have the following commuting diagram.
\[\begin{CD}
(G,k(\widehat{G'},G')) 		 @>>> 	(G,w)\\
@VVV       		 @VVV\\
(\widehat{G'},k(\widehat{G'},G'))	@>>>  (\widehat{G'},w(\widehat{G'},G')) \\
\end{CD}\]
where
the evaluation maps given by the vertical arrows are topological embeddings, and the identity maps given by the horizontal arrows are continuous.
By \ref{prelims}.8, the weak topology $w(\widehat{G'},G')$ and the locally compact topology 
$k(\widehat{G'},G')$ have the same collection of compact subsets on $\widehat{G'}$, which implies that $(G,k(\widehat{G'},G'))$ is itself hemicompact (since $(G,w)$ is hemicompact), and metrizable, as subgroup of $(\widehat{G'},k(\widehat{G'},G'))$. It follows then, by \cite{arens1946} again, that  $(G,k(\widehat{G'},G'))$ is a \lc \ subgroup of $(\widehat{G'},k(\widehat{G'},G'))$.
Therefore, we have proved that $(G,k(\widehat{G'},G'))$
is locally compact and metric. By Lemma \ref{lem_chacterization2}, its dual group is $G'$. By Theorem \ref{prelims}.1 it follows that $G=\widehat{G'}$ and 
since 
$(G,k(\widehat{G'},G'))\widehat{\:\:}=G'=(G,w)\widehat{\:\:}$ 
we have that $(G,k(\widehat{G'},G')^+)=(G,w)$, hence $(G,w) \in \sB$. 
This completes the proof. \epf}

In case $G$ has Ulam-measurable cardinality, Theorem \ref{th_characterization} does not hold because there are compact groups that are not $g$-sequential
(see \cite{ComfortRemus1994}). In this case, we must replace sequential continuity by $k$-continuity.


\bthm \label{th_characterization2}
Let $(G,w)$ be a hemicompact, totally bounded group.
Then $(G,w)\in \mathcal B$ \sii \ the following properties hold:
\begin{enumerate}
\item  Every $k$-continuous character on $(G,w)$ is continuous.
\item There exists a compact subgroup $K$ of $G$ such that $G/K$  has a countable network at $0$ consisting of $k$-neighborhoods of $0$.
\end{enumerate}
\ethm

\bexm\label{coun_th_characterization}
{The hemicompactness condition on Theorems  \ref{th_characterization} and \ref{th_characterization2} cannot be relaxed. If $G=\QQ^+$, then $G$ satisfies (1) and (2) but $G \not\in \sB$ since it is not hemicompact \cite{arens1946}.}
\eexm

In connection with this notion, the question of characterizing those totally bounded abelian groups $(G,w)$ such that $(G,kw)$ is locally compact and ($kw)^+ = w$ is proposed in \cite{comftrigwu}, and studied further by Galindo \cite{Galindo2004}. Next, we show how this question is related to the subject matter of this paper.

\blem \label{lem_Q4.3}
Let $(G,w)$ be a totally bounded group. Then $(G,w)\in \mathcal B$ \sii\ $(G,kw)$ is locally compact and ($kw)^+ = w$.
\elem

\pf\ Sufficiency is obvious. In order to prove necessity, assume that $(G,w)\in \mathcal B$.
That is, there is a locally compact topology $\tau$ on $G$ such that $\tau\supseteq w$ and $\tau^+=w$.
Since the topology of every locally compact group is both $g$-sequential and $k$-space
(see \cite{varop64}) and $\tau\supseteq w$, it follows that $\tau\supseteq s_g[w]\cup kw$.
On the other hand, by Theorem \ref{prelims}.8,
the groups $(G,w)$ and $(G,\tau)$ have the same collection of compact subsets,
which implies that $\tau\subseteq s_g[w]\cap kw$. In other words, we have that
$\tau=s_g[w]=kw$.
\epf

In light of the previous lemma, the next theorem provides an answer to Question 4.3 in \cite{comftrigwu}.

\bthm \label{th_characterization_general}
Let $(G,w)$ be totally bounded group whose cardinality is not Ulam-measurable.
Then $(G,w)\in \mathcal B$ \sii \ the following properties hold:
\begin{enumerate}
\item  Every sequentially continuous character on $(G,w)$ is continuous.
\item There is a hemicompact subgroup $(H,w)$ of $(G,w)$ \st \ $(G/H,w)= (G/H,\#)$
\item There exists a compact subgroup $K$ of $H$ such that $H/K$  has a countable network at $0$ consisting of $k$-neighborhoods of $0$.
\end{enumerate}
\ethm 


Again, when the group $G$ has Ulam-measurable cardinality, we have the following variant of Theorem \ref{th_characterization_general}.

\bthm \label{th_characterization_general2}
Let $(G,w)$ be totally bounded group.
Then $(G,w)\in \mathcal B$ \sii \ the following properties hold:
\begin{enumerate}
\item Every $k$-continuous character on $(G,w)$ is continuous.
\item There is a hemicompact subgroup $(H,w)$ of $(G,w)$ \st \ $(G/H,w)= (G/H,\#)$
\item There exists a compact subgroup $K$ of $H$ such that $H/K$  has a countable network at $0$ consisting of $k$-neighborhoods of $0$.
\end{enumerate}
\ethm
\mkp

We now establish the independence of the three conditions in Theorems \ref{th_characterization_general} and \ref{th_characterization_general2}.

\bexm\label{condicion1}
Set $G=(\ZZ,w)$, where $w$ is a totally bounded topology on $\ZZ$ such that $w\subsetneqq w(\ZZ,\TT)$ but contains no infinite compact subsets (see \cite{comftrigwu}).
Then $G$ satisfies (2,3) but fails to satisfy assertion (1). For, (1) does not hold since the only convergent sequences in $G$ are eventually constant; (2) $G$ is obviously hemicompact; and (3) $\{\{0\}\}$ is obviously a countable network at $\{0\}$ consisting of $k$-\nhd s.
\eexm

\bexm\label{condicion3}
Let $X$ be a compact metric space, $H_1=A(X)^+$, where $A(X)$ denotes the free Abelian group generated by $X$, and let $H_2=\RR^+$. Set $G=H_1\times H_2$. Then $G$ satisfies (1,2) but  fails to satisfy (3). For, a sequentially continuous character is obviously continuous on $X$, and thus it will be continuous by the properties of free Abelian groups. By \cite{Ark-Tka-book} (7.4.4) and (7.1.13) $A(X)$ is hemicompact, and by \cite{gal-her-fm1999} (4.20) $H_1$ is hemicompact as well; since $H_2$ is hemicompact by Theorem \ref{prelims}.8, a simple verification shows that $G$ is also hemicompact. Notice also that if $H_1$ satisfied (3),
then $\widehat{A(X)}$ would be first countable, hence \lc , which is absurd \cite{arens1946}. Thus, $G$ does not satisfy (3). 
\eexm

\bexm\label{condicion2}
If $G=\QQ^+$, then $G$ satisfies (1) and (3) but fails (2) \cite{arens1946}.
\eexm


\section{Respecting compactness}

Let us recall that a group is {\em von Neumann complete} if every closed precompact subset is complete.

\blem \label{Le_sequences}
Let $G$ be a  MAP (von Neumann) complete group that respects compactness. If $(x_n)$ is a Cauchy sequence in $G^+$,
then it converges in $G$.
\elem 

\pf\ If $(x_n)$ is precompact in $G$, then $\overline{(x_n)}^G$ is a compact subset in $G$ homeomorphic
to $\overline{(x_n)}^{bG}$. As a consequence $(x_n)$ must be convergent in $G$. Therefore, we may assume
that $(x_n)$ is not precompact. Furthermore, taking a convenient subsequence if it were necessary,
we may assume that $(x_n)$ is uniformly discrete. That is, there is a \nhd \ of the identity, say $U$, such that
$x_nx^{-1}_m\notin U$ for all $n,m$ in $\NN$. Now, since $(x_n)$ is Bohr Cauchy, it follows that
$(x_nx^{-1}_{n+1})$ Bohr converges to the neutral element. Since $G$ respects compactness, we have that
$(x_nx^{-1}_{n+1})$ converges to the neutral element in $G$, which is a contradiction, completing the proof.\epf

As a consequence of the previous lemma, we obtain:

\blem \label{Le_CTW}
Let $G$ be a MAP (von Neumann) complete group that respects compactness, and let $N$ be a closed, metrizable subgroup of $bG$.
Set $H^+=G^+\cap N$. Then $H$, the inverse image of $H^+$ in $G$, is a compact metrizable group isomorphic to $H^+$.
\elem 


\begin{defn} \label{str_resp_cpt} { We say that a MAP group  $G$ {\em strongly respects compactness} if whenever $N$ is a closed {\em metrizable} subgroup of the Bohr compactification $bG$ of $G$ and $A \subseteq G$, then $A + (N \cap G)$ is compact in $G$ whenever $\phi(A)$ is compact,  where $\phi$ is the obvious map $G \to bG \to bG/N$. If $N$ is a closed  subgroup of $bG$ \st \ for any $A \subseteq G$, $A + (N \cap G)$ is compact in $G$ whenever $\phi(A)$ is compact in $bG/N$, then we say that $N$ {\em preserves compactness}.}\end{defn}

\bthm \label{Th_Polish}
 Let $G$ be a MAP Abelian Polish group (more generally, a MAP abelian metrizable von Neumann complete group)
 that respects compactness. Then $G$ strongly respects compactness.
\ethm 

\pf\ Let $N$ be closed, metrizable subgroup of $bG$ and assume that $A \subseteq G$ with $\phi(A)$ compact in $bG/N$.
We have to prove that $A+(N\cap G)$ is compact in $G$. Now, since $A+(N\cap G)$ closed
in $G^+$, it is also closed in $G$, which is complete. Therefore,
it will suffice to show that $A+(N\cap G)$ is precompact in $G$. Assume 
{ otherwise}.
Then $A+(N\cap G)$ must contain an infinite uniformly sequence $(x_n)=(a_n+y_n)$, where
$(a_n)\subseteq A$ and $(y_n)\subseteq N\cap G$. 
Now, $A+N$ is a compact metrizable subspace of $bG$ (in order to see this, use the first part of the proof of Lemma 2.6 in \cite{comftrigwu}).
Therefore, $(x_n)$ must contain a Cauchy subsequence, say $(x_{n_m})$. By Lemma \ref{Le_sequences}, this subsequence
converges to a point $p\in G$, which contradicts our assumption about $(x_n)$ being uniformly discrete.
Thus $A+(N\cap G)$ is precompact in $G$, which completes the proof.\epf

{The following is Question 4.1 in \cite{comftrigwu}: 

\begin{quest} \label{CTW4.1} Let $G$ be a MAP group and suppose that every closed, metrizable
subgroup $N$ of $bG$ such that $N \cap G = \{0\}$ preserves compactness.
Does it follow that $G$ strongly respects compactness?
\end{quest} }

We will need the following:

\begin{defn} \label{k_omega}
{ A Hausdorff topological space $X$ is a {\em $k_\omega$-space} if there exists an ascending sequence of compact subsets $K_1 \subseteq K_2\subseteq \cdots X$ such that $X = \cup_{n<\omega} K_n$ and
$U \subseteq X$ is open if and only if $U \cap K_n$ is open in $K_n$ for each $n < \omega$.}
\end{defn}

The following result answers { Question \ref{CTW4.1} in the negative.}

\bexm \label{Ex_q1}
{ Take $G_n:=\TT$ and $H_n:=\QQ/\ZZ$ for all $n<\omega$, and set $G:=\bigoplus\limits_{n<\omega} G_n$,
and $H:=\bigoplus\limits_{n<\omega} H_n$, with both groups
equipped with the 
box topology}. From here on, we identify the groups $G_n$, $H_n$,
$\bigoplus\limits_{n<N} G_n$ and $\bigoplus\limits_{n<N} H_n$
with their isomorphic subgroups in $G$ and $H$ respectively.
We have the following facts:\eexm

\begin{enumerate}
\item $G$ { is the countable direct limit of compact groups and, therefore it is a MAP, $k_\omega$-group having
$H$ as a dense subgroup.} 
\item $G$ strongly respects compactness.
\item For every compact subset $K$ of $G$ there is $n_0\in\NN$ such that  $K\subseteq \bigoplus\limits_{n<n_0} G_n$.
\item $bG=bH$.
\item If $N$ is a compact {metrizable} subgroup of $bG$ such that $N\cap H=\{0\}$, then $N\cap G=\{0\}$.
\item If $N$ is a compact {metrizable} subgroup of $bG$ such that $N\cap H=\{0\}$, then $N$ respects compactness {in $H$}.
\item$H$ does not strongly respects compactness.
\end{enumerate}

\pf\ (1) is clear. (2) 
is proved in \cite{FHT:iii},
where it is established that every locally $k_\omega$-group strongly respects compactness.
(3) is clear since $G$ is equipped with the countable box topology. (4) follows from the density of
$H$ in $G$. In order to prove (5), reasoning by contradiction, assume that $N\cap G\not=\{0\}$. Since $G$ strongly respects compactness,
it follows that $N\cap G$ is a compact subgroup of $G$ and, by (3), there is $n_0\in\NN$ such that
$G\cap N\subseteq \bigoplus\limits_{n<n_0} G_n$. Now, every proper closed subgroup of $\TT$ is finite and contained
in $\QQ/\ZZ$. Thus, if $\pi_n$ denotes the $nth$-projection {\bf of $G\cap N$} onto $G_n$, we have that either
$\pi_n(G\cap N)=\TT$ or $\pi_n(G\cap N)\subseteq \QQ/\ZZ$. If Ker\,$\pi_1$ were trivial, then $\pi_1$ is $1$-to-$1$ and, therefore,
$\pi_1(G\cap N)$ contains a finite subgroup $F$ that is isomorphic to its inverse image
$\pi_1^{-1}(F)\subseteq H$. Thus, we may assume that Ker\,$\pi_1$ is nontrivial. Then we replace $G\cap N$ by Ker\,$\pi_1$.
Applying induction, it follows that $H\cap N\not=\{0\}$. (6)
Let $N$ be a compact metrizable subgroup of $bH$ such that $H\cap N=\{0\}$ and let $A\subseteq H$
such that $A+N$ is compact in $bH$. By (2) and (4) we have that $A+(G\cap N)$ is compact in $G$ and, by (5),
$G\cap N=\{0\}$. Thus $A=A+(G\cap N)$ is compact in $G$. Since $A\subseteq H$, we obtain that $A$
is compact in $H$. In other words, the group $N$ respects compactness in $H$ if $H\cap N=\{0\}$.
(7) Take $N=G_1$ and $A=H_1$. If $\phi\colon bH\to bH/N$ denotes the canonical quotient map,
we have that $\phi(A)=\{0\}$ is trivially compact in $bH/N$. On the other hand, $A+(N\cap H)=H_1$,
which is not compact in $H$. Therefore $H$ does not strongly respects compactness.\epf


\bibliographystyle{amsplain}
\bibliography{weakbohr}

\providecommand{\bysame}{\leavevmode\hbox to3em{\hrulefill}\thinspace}
\providecommand{\MR}{\relax\ifhmode\unskip\space\fi MR }
\providecommand{\MRhref}[2]{%
  \href{http://www.ams.org/mathscinet-getitem?mr=#1}{#2}
}
\providecommand{\href}[2]{#2}
\begin{thebibliography}{10}

\bibitem{arens1946}
R.~Arens, \emph{A topology for spaces of transformations}, Ann. Math.
  \textbf{47} (1946), 480--495.

\bibitem{Ark-Tka-book}
{A}. {A}rhangel'skii and {M}. {T}kachenko, \emph{Topological {G}roups and
  {R}elated {S}tructures}, vol.~1, 2008, Atlantic Press, Amsterdam-Paris.

\bibitem{comf-hbook}
W.~W. Comfort, \emph{Topological groups}, Handbook of {S}et-{T}heoretic
  {T}opology (K.~Kunen and J.~E. Vaughan, eds.), Elsevier Science Publishers,
  B. V., Amsterdam, 1984, pp.~1143--1263.

\bibitem{ComfortRemus1994}
W.~W. Comfort and D.~Remus, \emph{Compact groups of {U}lam-measurable
  cardinality: partial converses to theorems of {A}rhangel'ski\u\i \ and
  {V}aropoulos}, Math. Japon. \textbf{39} (1994), no.~2, 203--210.

\bibitem{ComfortRoss1964}
W.~W. Comfort and K.~A. Ross, \emph{Topologies induced by groups of
  characters}, Fundamenta Mathematicae \textbf{55} (1964), no.~3, 283--291.

\bibitem{comftrigwu}
W.~W. Comfort, F.~J. Trigos-Arrieta, and T.-S. Wu, \emph{{The Bohr
  compactification , modulo a metrizable subgroup}}, Fund. Math. \textbf{143}
  (1993), 119--136.

\bibitem{FHT:iii}
M.~V. Ferrer, S.~Hern{\'{a}}ndez, and L.~T{\'{a}}rrega, \emph{{Interpolation
  sets in spaces of continuous metric-valued functions}}, J. Math. Analysis and
  Applications \textbf{466} (2018), no.~2, 1426--–1442.

\bibitem{Galindo2004}
J.~{Galindo}, \emph{Totally bounded group topologies that are {B}ohr topologies
  of {LCA} groups}, {T}opology {P}roceedings \textbf{28(2)} (2004), 467--478.

\bibitem{gal-her-fm1999}
J.~{Galindo} and S.~{Hern\'andez}, \emph{The concept of boundedness and the
  {B}ohr compactification of a {MAP} {A}belian group}, Fund. Math. \textbf{159}
  (1999), 195--218.

\bibitem{GalHerWu}
J.~{Galindo}, S.~{Hern\'andez}, and Ta-Sun {Wu}, \emph{Recent results and open
  questions relating {C}hu duality and {B}ohr compactifications of locally
  compact groups}, Open problems in topology. II. Edited by: Elliott Pearl.
  Amsterdam: Elsevier (2007).

\bibitem{glicks1962}
I.~Glicksberg, \emph{Uniform boundedness for groups}, Canadian J. Math.
  \textbf{14} (1962), 269--276.

\bibitem{Hernandez1998}
S.~Hern\'andez, \emph{The {D}imension of an {LCA} group in its {B}ohr
  {T}opology}, Topology and its Applications \textbf{86} (1998), 63--67.

\bibitem{HT_jmaa_2005}
S.~Hern\'{a}ndez and F.~J. Trigos-{A}rrieta, \emph{Group duality with the
  topology of precompact convergence}, J. Math. Analysis and Applications
  \textbf{303} (2005), 274--287.

\bibitem{hewitt1963abstract}
E.~Hewitt and K.~A. Ross, \emph{{A}bstract {H}armonic {A}nalysis: Volume {I},
  {S}tructure of {T}opological {G}roups, {I}ntegration {T}heory, {G}roup
  {R}epresentations}, Grundlehren der mathematischen Wissenschaften, Springer
  New York, 1963.

\bibitem{Heyer1970}
H.~Heyer, \emph{{Dualit{\"{a}}t lokalkompakter Gruppen}}, Springer-Verlag,
  Berlin-Heidelberg-New York, 1970.

\bibitem{Hughes1973}
R~Hughes, \emph{{Compactness in locally compact groups}}, Bulletin of the
  American Mathematical Society \textbf{79} (1973), 122--123.

\bibitem{noble70b}
N.~Noble, \emph{$k$-groups and duality}, Trans. Amer. Math. Soc. \textbf{151}
  (1970), 551--561.

\bibitem{racz-trig2001}
S.~U. Raczkowski and F.~J. Trigos-Arrieta, \emph{Duality of {T}otally {B}ounded
  {G}roups}, Bolet\'\i n de la Sociedad Matem\'atica Mexicana \textbf{7}
  (2001), 1--12.

\bibitem{reid}
G.~A. Reid, \emph{On sequential convergence in groups}, Math. Zeitschrift
  \textbf{102} (1967), 225--235.

\bibitem{Trigos-Arrieta1991}
{F}rancisco~{J}avier Trigos-Arrieta, \emph{{Pseudocompactness on groups}}, Ph.D
  Thesis, Wesleyan University (1991).

\bibitem{varop64}
N.~Th. Varopoulos, \emph{A theorem on the continuity of homomorphisms of
  locally compact groups}, Proc. Cambridge Phil. Soc. \textbf{60} (1964),
  449--463.

\end{thebibliography}
\end{document}